\parskip=12pt plus 6pt
\tolerance=2500

\magnification=\magstep1

\input psfig.sty

\font\eightbf=cmbx8
\font\bigbf=cmbx10 scaled \magstep1

\font\vbigbf=cmbx10 scaled \magstep2

\font\tenmsb=msbm10
\font\sevenmsb=msbm7
\font\fivemsb=msbm5
\newfam\msbfam
\textfont\msbfam=\tenmsb
\scriptfont\msbfam=\sevenmsb
\scriptscriptfont\msbfam=\fivemsb
\def\Bbb{\fam\msbfam}

\def\ZZ{{\Bbb Z}}

\def\RR{{\Bbb R}}

\def\SS{{\Bbb S}}

\def\itembu{\par\noindent{$\bullet$}\hskip10pt}
\def\nobf{\noindent\bf}
\def\nobigbf{\noindent\bigbf}

\def\emptyset{/\kern -0.67em{\scriptstyle\bigcirc}}

\def\bbuildrel#1_#2^#3{\mathrel{
\mathop{\kern 0pt#1}\limits_{#2}^{#3}}}

\centerline{\vbigbf An overview of domino and lozenge tilings}
\smallskip
\centerline{\it Nicolau C. Saldanha and Carlos Tomei}
\bigskip

{

\narrower

{\nobf Abstract:}
We consider tilings of quadriculated regions by dominoes
and of triangulated regions by lozenges.
We present an overview of results concerning tileability,
enumeration and the structure of the space of tilings.

}

\goodbreak
\bigskip
{\nobigbf Introduction}

\nobreak
\smallskip
\nobreak

Let $R$ be a finite juxtaposition of unit squares in the plane.
We consider tilings of $R$ by {\sl dominoes},
$2 \times 1$ rectangles.
Alternatively, if $R$ is obtained from unit equilateral triangles,
we try to tile $R$ by {\sl lozenges} (or {\sl calissons}),
unions of two adjacent triangles.
The following questions are natural:
{\parskip=6pt
\itembu{Can $R$ be tiled?}
\itembu{How many tilings exist?}}

A next set of questions is suggested by the concept of {\sl flip}.
We perform a flip by lifting two dominoes forming a square
or three lozenges forming a hexagon
and placing them back after a rotation of $90^\circ$ or $60^\circ$.
{\parskip=6pt
\itembu{Do all tilings admit flips?}
\itembu{Are any two tilings joined by a sequence of flips?}
\itembu{What is the distance (the minimal number of flips)
between two tilings?}
\itembu{Which flips decrease the distance to a given tiling?}
\itembu{Are there closed geodesics of flips in the space of tilings?}}

Furthermore, as we shall see, there is a natural bipartition
on the space of domino tilings so that flips
connect tilings in opposite classes.
{\parskip=6pt
\itembu{Is the number of tilings in each class the same?}}

For simply connected regions, these questions are easier.
At the other end of the spectrum, they still make sense
for quadriculated or triangulated surfaces.

Section 1 deals with the problem of tileability;
the main references are [CL] and [T].
In section 2, we discuss some techniques for counting tilings,
from the beginning of the century with [M]
to recent contributions by mathematical physicists ([K] and [LL]).
Finally, section 3 deals with flips and related issues;
the results were originally presented in [STCR].

\goodbreak
\bigskip
{\nobigbf 1. Existence of tilings}

\nobreak\smallskip\nobreak

Given a finite set of types of tiles, can we use them to tile the plane?
Rather surprisingly, this is an undecidable question ([B]).
The corresponding problem for finite regions is NP-complete ([GJ]).
Colouring arguments are among the simplest
to provide necessary conditions for tileability.
Both the triangulated and the quadriculated plane are bicolourable
and a region must have the same number of triangles
(or squares) of each colour to be tiled by lozenges (or dominoes).
Conway and Lagarias ([CL]) discovered a method to obtain
strong necessary conditions for tile\-ability.
Following Thurston ([T]), we present their method
for simply connected quadriculated or triangulated regions,
the tiles being dominoes or lozenges.

\midinsert
{
\smallskip

\centerline{\psfig{file=fig11.eps,width=2.5in}}


\smallskip
\centerline{\eightbf Figure 1.1}

}
\endinsert

We start with lozenges.
In the triangulated plane, we may label the oriented edges
$a$, $b$ and $c$ as in figure 1.1.a and paths along edges yield words
in the letters $a$, $b$ and $c$.
Notice that a path is closed iff the total numbers of $a$'s,
$b$'s and $c$'s are equal, i.e., iff its word equals the identity
in $G_1 = \ZZ^2$, given in terms of generators and relations as
$G_1 = \langle a, b, c | abc = acb = e \rangle$.
Consider the group
$G = \ZZ^3 = \langle a, b, c | [a,b] = [a,c] = [b,c] = e \rangle$
(as usual, $[g,h] = g^{-1}h^{-1}gh$):
the relations are the boundaries of the three possible
orientations of lozenges (i.e., position up to translation).
The boundary of a simply connected region $R$ yields a word $w$
and thus an element $g \in G$:
a tiling of $G$ gives a way to write $w$ as a juxtaposition
of conjugates of boundaries of tiles
(see figure 1.2 or [CL] for an actual proof)
and thus $g = e$.

\midinsert
{
\smallskip

\centerline{\psfig{file=fig12.eps,width=4.5in}}


\smallskip
\centerline{\eightbf Figure 1.2}

}
\endinsert

These group-theoretical considerations produce
Conway-Lagarias's condition for tile\-ability
of simply connected regions.
We state their result for triangulated (or quadriculated) regions:

{\nobf Theorem 1.1: }
{\sl Let $R$ be a triangulated (or quadriculated)
simply connected bounded region and
$\Sigma$ a finite set of (connected and simply connected) tile shapes.
Let $G$ be the group with generators corresponding to edges
and relations given by boundaries of elements of $\Sigma$.
Then a necessary condition for the tileability of $R$
by translates of elements of $\Sigma$ is that the word
induced by the boundary of $R$ is trivial in $G$.}

For lozenges, this condition reduces
to the simple colour condition.
More interesting applications of the group theoretic conditions
yield the following results ([CL]).
Let $T_n$ be the triangular region of size $n$ in the hexagonal
lattice (fig 1.3.a) and $L_3$ be the {\sl tribone},
the juxtaposition of three hexagons in a line (fig 1.3.b).

{\nobf Theorem 1.2: }
{\sl The triangular region $T_n$ can be tiled
by copies of $T_2$ iff $n \equiv 0, 2, 9$ or $11 \pmod{12}$.
Also, $T_n$ cannot be tiled by congruent copies of $L_3$.}

\midinsert
{
\smallskip

\centerline{\psfig{file=fig13.eps,width=3.5in}}


\smallskip
\centerline{\eightbf Figure 1.3}

}
\endinsert

David and Tomei noticed that a lozenge tiling (as in figure 1.4.a)
may literally be seen as a picture of a bunch of cubes
stored in a big box (as in figure 1.4.b).
Align the sides of the box with the coordinate axes
so that vertices of triangles receive integer $x, y, z$ coordinates.
In figure 1.4, for example, the box
has size $5 \times 5 \times 5$
and there are 66 cubes in it.
The figure is an orthogonal projection
of the three dimensional configuration on the plane $x+y+z = 0$.
If we think of gravity as pointing ``downwards'' to $(-1,-1,-1)$,
the cubes form a gravitationally stable {\sl pile}
iff they represent a tiling.
By projecting on the coordinate planes,
we immediately obtain for hexagonal regions
the following invariance result ([DaT]):

\midinsert
{
\smallskip

\centerline{\psfig{file=fig14a.eps,width=2.0in}\hglue0.5in%
\psfig{file=fig14b.eps,width=2.0in}}


\smallskip
\centerline{\eightbf Figure 1.4}

}
\endinsert

{\nobf Theorem 1.3: }
{\sl Let $R$ be an arbitrary triangulated planar region.
Then the numbers $n_{xy}, n_{xz}, n_{yz}$ of
lozenges in each orientation do not depend on the tiling.}

The visual argument works for simply connected regions;
for general regions we might use {\sl height sections}
(to be introduced later) but we present an elementary proof
which appears to be folklore.
Consider vectors going from the origin of the plane
to the centre of each white triangle
and from the centre of each black triangle
to the origin and add them up obtaining a vector $v$.
Let $w_{xy}$ be the vector joining the centres of triangles
in a $xy$-lozenge, from black to white;
define $w_{xz}$ and $w_{yz}$ similarly.
Given a tiling, we may group pairs of vectors
associated to the same lozenge to obtain:
$$v = n_{xy}w_{xy} + n_{xz}w_{xz} + n_{yz}w_{yz}.$$
Since $n_{xy} + n_{xz} + n_{yz}$ does not depend on the tiling
and $w_{xy} - w_{xz}$ and $w_{xy} - w_{yz}$ are 
linearly independent, the theorem follows.

Thurston obtained the same three dimensional interpretation
as a spinoff of the group-theoretical construction.
Consider the Cayley diagram $\Gamma$ for $G$ and its natural projection
onto $\Gamma_1$, the Cayley diagram of its quotient $G_1$
($\Gamma_1$ is, of course, the triangulated plane):
a tiling of a simply connected region lives in $\Gamma_1$
but may be lifted to $\Gamma$ (since each tile can).
Imbed the Cayley graph $\Gamma$ of $G=\ZZ^3$ in $\RR^3$
(in the obvious way) and $\Gamma_1$ in the plane $x+y+z=0$:
the projection of $\Gamma$ onto $\Gamma_1$
is the restriction of the orthogonal projection onto this plane. 
This lifting assigns to each vertex of a triangle
the same coordinates as the more visual procedure
represented in figure 1.4.

\midinsert
{
\smallskip

\centerline{\psfig{file=fig15a.eps,width=3.0in}\hglue1in
\raise1.1in\hbox{\psfig{file=fig15b.eps,width=1.2in}}}


\smallskip
\centerline{\eightbf Figure 1.5}

}
\endinsert

The extra information obtained by lifting is $x+y+z$,
the height with respect to the plane $x+y+z=0$:
we call the function which assigns to each vertex
the sum of its coordinates a {\sl height function},
as shown in figure 1.5.a.
Of course, there are other expressions
which convey the 3D information equally well
but, for height functions of lozenge tilings, we always take $x+y+z$.
Height functions are therefore well defined
for a tiling by lozenges of a simply connected region
up to an additive constant.
For a given tiling, a simple local rule obtains the height function:
when following an edge at the border of a tile,
add (resp. subtract) 1 to the height function at the source
to get the height function at the target
if the edge has a black triangle to the right (resp. left).
Also, height functions of a simply connected region $R$
are easy to characterize:
height functions have the mod 3 values shown in figure 1.5.b,
change by at most 2 along any edge in $R$
and change by 1 along an edge in the boundary of $R$.
An unexpected consequence of this characterization
is the following ([STCR]):

{\nobf Theorem 1.4: }
{\sl
For a simply connected region,
let $T$ be the space of its lozenge tilings.
The set $T$ is a lattice:
given two tilings, the maximum (or minimum) of their height
functions is the height function of some tiling.}

In terms of the pile-of-cubes interpretation,
the minimum tiling corresponds to the empty box.
Thurston ([T]) uses height functions to prove:

{\nobf Theorem 1.5: }
{\sl There is a polynomial time algorithm
to decide whether a simply connected region $R$
is tileable by lozenges.}

In a nutshell,
start by computing the height function along the boundary of $R$
and move inwards assigning to each point the smallest value
satisfying the characterization.
The algorithm either fills the region, obtaining the minimal tiling,
or runs into a contradiction, indicating non-tileability of $R$.

An equivalent criterion, computationally harder
but easier to state, is the following:
consider positively oriented simple closed curves
which, outside the boundary of the region,
have a black triangle to their left---%
the numbers of enclosed black and white triangles are equal.

The group theoretical construction easily applies to domino tilings.
Label oriented edges in the quadriculated plane as in figure 1.1.b,
obtaining the Cayley diagram $\Gamma_1$ for
$G_1 = \langle a, b | [a,b] = e \rangle = \ZZ^2$.
Taking into account the shape of the tiles, define
$G = \langle a, b | [a,b^2] = [a^2,b] = e \rangle$:
it turns out that $G$ is a non-abelian group and we describe
a construction of its Cayley diagram $\Gamma$ in $\RR^3$,
obtaining a height function for domino tilings.

Introducing a third generator $c = [a,b]$, we have
$G = \langle a, b, c | ac = c^{-1}a, bc = c^{-1}b, ab = bac \rangle$
and, by standard arguments, any element of $G$
can be written as $c^w b^y a^x$ for unique integers $x, y, w$.
We may thus take $\ZZ^3$ to be the set of vertices
of a Cayley graph for $G$:
connect $(x,y,w)$ by an $a$-edge to $(x+1,y,w)$
(because $c^wb^ya^x \cdot a = c^wb^ya^{x+1}$),
by a $b$-edge to $(x,y+1,w)$ if $x$ is even or
to $(x,y+1,w+(-1)^{y})$ if $x$ is odd
and by a $c$-edge to $(x,y,w+(-1)^{x+y})$.
A somewhat nicer picture is obtained by taking
$$z = 4w + \cases{
0 & if $x$ and $y$ are both even,\cr
1 & if $x$ is odd and $y$ is even,\cr
-2 & if $x$ and $y$ are both odd,\cr
-1 & if $x$ is even and $y$ is odd.\cr}$$
The vertices of $\Gamma$ are then the points $(x,y,z) \in \ZZ^3$
for which
$$z {\bbuildrel{\equiv}_{\pmod 4}^{}}
\cases{
0 & if $x$ and $y$ are both even,\cr
1 & if $x$ is odd and $y$ is even,\cr
2 & if $x$ and $y$ are both odd,\cr
3 & if $x$ is even and $y$ is odd,\cr}$$
and, for the original presentation of $G$,
there are edges joining two vertices
iff the euclidean distance is $\sqrt{2}$.
Notice that this graph projects onto $\Gamma_1$
by omitting the $z$ coordinate.

By construction, any domino tiling lifts to this Cayley graph.
In figure 1.4 we show a tiling and its related height function.
The simple procedure to obtain height functions is:
colour squares of the region in a checkerboard pattern,
assign 0 to the origin and, for an edge not covered by a domino,
assign to the target vertex the value
at the source plus (resp. minus) 1
if the square to the left is black (resp. white).

\midinsert
{
\smallskip

\centerline{\psfig{file=fig16.eps,width=2.0in}}


\smallskip
\centerline{\eightbf Figure 1.6}

}
\endinsert

Again, for simply connected regions,
it is easy to give a local characterization of height functions,
from which the counterparts of theorems 1.4 and 1.5 follow 
([STCR] and [T]).

The pile-of-cubes interpretation for lozenge tilings
has its analogue for domino tilings.
The vertices of the three dimensional solid involved
belong to $\Gamma$ as above;
the solid itself can be taken to be
$$\left\{{(x,y,z)\;\big|\;
|x| \le 1, |y| \le 1,
|x| + 3 |y| - 2 |xy| \le z \le 4 - 3 |x| - |y| + 2 |xy|}\right\},$$
shown in figure 1.7.
Notice that, when seen from the top (resp. bottom),
the solid covers two vertical (resp. horizontal) dominoes.
We call this solid a block and there is, indeed,
a {\sl pile-of-blocks interpretation};
height functions, however, nicely encapsulate
the relevant properties of such arrangements 
and we may therefore avoid visualization altogether.

\midinsert
{
\smallskip

\centerline{\psfig{file=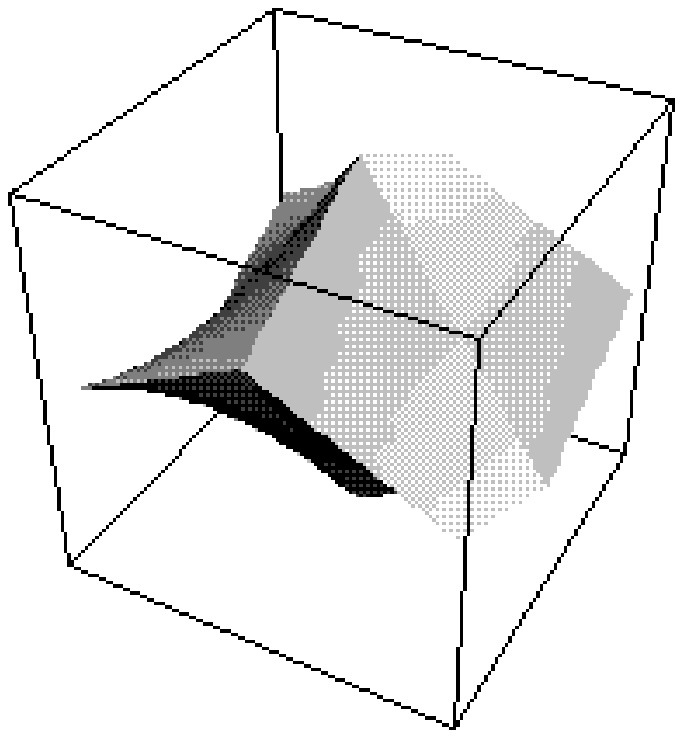,width=3.5in}}


\smallskip
\centerline{\eightbf Figure 1.7}

}
\endinsert

The procedures in [CL] and [T] do not apply to regions
which are not simply connected.
To study domino tilings,
we introduced {\sl height sections} ([STCR]),
a generalization of height functions which applies
to any quadriculated surface (e.g., a torus or Klein bottle).
Instead of a formal definition, we give a typical example;
vertices in the interior of the surface
must be surrounded by 4 squares.
A similar construction applies for lozenge tilings
of triangulated surfaces
(interior vertices in a triangulated surface
are surrounded by 6 triangles).

Consider a cylinder with boundary, as in figure 1.8, tiled by dominoes.
The two vertical sides are identified as indicated by the labels:
two dominoes trespass the cut.
Begin by assigning the value 0 to the vertex on the left labelled $A$
and follow the local instructions for constructing height functions
to obtain the value of the height section at all vertices.
The point $A$ in the cylinder is represented by two vertices
in the diagram which received different numbers:
0 at the left and 2 at the right.
The same ambiguity happens for the other points along the cut.
How is this to be interpreted?
We might take values mod 2 but this would amount
to throwing away the bundle with the section.
Instead, we construct a fibre bundle with fibre $\ZZ$
and the cylinder as base space (a 2-manifold with boundary).
Over any simply connected open subset of the basis,
the leaves can be labelled by the integers:
do this for the complement of the cut---%
the values of the height section are to be understood
in terms of these local coordinates.
The different values for the points in the cut indicate
how to glue leaves along the cut:
in our example, the $n$-leaf at the left
is glued to $n+2$-leaf at the right.
Each number in the figure labels one element
of the fibre over the corresponding vertex 
and the different numbers over distinct representations
of the same vertex label, by construction, the same element.
These elements form a section of the restriction
of this {\sl height bundle} to the vertices of the basis:
the {\sl height section}.
Summing up, the procedure simultaneously constructs
the height section and the height bundle.

\midinsert
{
\smallskip

\centerline{\psfig{file=fig18.eps,width=2.0in}}


\smallskip
\centerline{\eightbf Figure 1.8}

}
\endinsert

Height sections allow us
to give a tileability algorithm similar to Thurston's.
It would be interesting to detail such an efficient algorithm
and generalize the alternate criterion (based on closed curves)
to arbitrary quadriculated surfaces.
For planar regions, the space of tilings is still a lattice.

\goodbreak
\bigskip
{\nobigbf 2. Counting tilings}

\nobreak\smallskip\nobreak

For a quadriculated or triangulated planar region $R$,
consider a graph $\Lambda$ whose vertices
are the $k$ unit faces (squares or triangles) of $R$,
connected by an edge iff the faces have a common side.
Since faces in $R$ can be coloured black and white,
$\Lambda$ is bipartite;
we assume the number of black faces to equal the number
of white faces, a necessary condition for the existence of tilings.
Let $A$ be the $k \times k$ adjacency matrix of $\Lambda$:
by labelling black vertices first, we may write
$$A = \pmatrix{0 & B \cr {B^t} & 0 \cr}.$$
In a language familiar to physicists,
tilings of $R$ correspond to {\sl dimer coverings} of $\Lambda$,
i.e., partitions of the set of vertices
into pairs of adjacent vertices.

Each dimer covering corresponds to a monomial
in the expansion of $\det B$.
Suppose that $\tilde B$ is obtained from $B$
by changing the signs of some entries and that
all monomials in the expansion of $\det \tilde B$ have the same sign.
Then the number $N$ of tilings of $R$ is $|\det \tilde B\,|$.
Amazingly, for simply connected triangulated regions,
we can take $\tilde B = B$.
For quadriculated regions,
Kasteleyn ([K]) showed that such a $\tilde B$ can be
obtained from $B$ by
changing the signs of the 1's corresponding to edges
on alternate vertical lines.
Lieb and Loss ([LL]) generalized this idea to prove
a formula for the number of dimer partitions
of an arbitrary bipartite planar graph.
The following theorem, a special case of their results,
unifies the previous claims.

{\nobf Theorem 2.1: }
{\sl
Let $\Lambda$ be a planar bipartite graph with adjacency matrix
$$A = \pmatrix{0 & B \cr {B^t} & 0 \cr}.$$
Choose a set of edges in $\Lambda$ such that
on each minimal loop of size $4n$ (resp. $4n+2$)
there is an odd (resp. even) number of chosen edges.
Change the sign of the 1's in $B$
corresponding to chosen edges obtaining $\tilde B$.
Then the number of dimer partitions of $\Lambda$ is $|\det \tilde B\,|$.
}

The required choice of edges is always possible.
An example is given in figure 2.1,
where there are 6 minimal loops of sizes 8, 4, 4, 4, 4 and 10
and the chosen edges are thick.

\midinsert
{
\smallskip

\centerline{\psfig{file=fig21.eps,width=3.5in}}


\smallskip
\centerline{\eightbf Figure 2.1}

}
\endinsert

The sign of the monomial in $\det B$ associated to a tiling
is not natural, depending on the labelling of vertices.
Given two tilings, however, equality of their signs is
well defined and thus tilings naturally split in two classes.
This splitting is very different
in the quadriculated and triangulated cases:

{\nobf Theorem 2.2: }
{\sl Let $R$ be a simply connected region,
$N$ the number of tilings of $R$ and $B$ as above.
If $R$ is quadriculated,
$$|\det B\,| \le 1:$$
the two classes are almost balanced.
If $R$ is triangulated,
$$|\det B\,| = N:$$
one of the classes is empty.
}

The quadriculated case is due to Deift and Tomei ([DeT]).
For general planar regions this difference assumes any value.
The triangulated case is a corollary of theorem 2.1
but we will see another proof in the next section.

We present closed formulae for $N$ in special cases.
Kasteleyn ([K]) obtained a closed form for the number of tilings
of a rectangle:

{\nobf Theorem 2.3: }
{\sl If $R$ is a $n \times m$ rectangle,
$$\prod_{k=1}^{m/2}\prod_{\ell=1}^{n}
{2\left({\cos^2{k\pi \over m+1} +
\cos^2{\ell\pi \over n+1}}\right)^{1/2}},$$
assuming $m$ even.}

Kasteleyn's proof boils down
to computing the determinant of $\tilde B$, as above.
He essentially finds the eigenvalues and eigenvectors of 
$$\pmatrix{0 & \tilde B \cr \tilde B^t & 0 \cr}^2;$$
the eigenvalues are the terms in the product
with doubles multiplicities
and the eigenvectors have the form
$$\sin{k \pi i \over m+1}\sin{\ell \pi j \over n+1}.$$

The pile-of-cubes interpretation of lozenge tilings
reduces their counting for center-symmetric hexagonal regions
to that of {\sl plane partitions}.
To describe a tiling, count the number of cubes above
each square in the $xy$ plane
(here, the ``up'' direction is $(0,0,1)$ and not $(1,1,1)$).
The example in figure 1.3 would be encoded as below.
$$\matrix{
5 & 5 & 5 & 5 & 3\cr
5 & 5 & 5 & 5 & 3\cr
5 & 5 & 1 & 1 & 0\cr
4 & 1 & 0 & 0 & 0\cr
3 & 0 & 0 & 0 & 0\cr
}$$
Following Bender and Knuth ([BK]),
a plane partition of $n$ (the total number of cubes in the pile)
is an array of non-negative integers
$$\matrix{
n_{11} & n_{12} & n_{13} & \cdots \cr
n_{21} & n_{22} & n_{23} & \cdots \cr
\vdots & \vdots & \vdots & \ddots \cr}$$
for which $\sum_{ij}n_{ij} = n$,
$n_{ij} \ge n_{(i+1)j}$ and $n_{ij} \ge n_{i(j+1)}$
(these conditions correspond to the pile being gravitationally stable).
If $n_{ij} = 0$ for all $i > r$ (resp. $j > c$),
the partition is called {\sl $r$-rowed} (resp. {\sl $c$-columned})
(this is the size of the region in the $x$ (resp. $y$) direction).
If $n_{11} \le m$, we say the parts {\sl do not exceed $m$}
(the size in the $z$ direction).
Let $\pi_n$ be the number of plane $r$-rowed,
$c$-columned partitions of $n$ with parts not exceeding $m$.
MacMahon ([M]) obtained
$$\sum_{n=0}^{\infty}{\pi_n x^n} =
\prod_{i=1}^r\prod_{j=1}^m{1 - x^{c+i+j-1} \over 1 - x^{i+j-1}},$$
yielding the following formula:

{\nobf Theorem 2.4: }
{\sl
The number of lozenge tilings
of a center-symmetric hexagonal region of sides $r$, $c$ and $m$ is
$$\prod_{i=1}^r\prod_{j=1}^m{c+i+j-1 \over i+j-1}.$$}

\goodbreak
\bigskip
{\nobigbf 3. Flips and the space of tilings}

\nobreak\smallskip\nobreak

The space $T$ of all tilings of a given region $R$
admits a natural graph structure:
two tilings are adjacent iff they differ by a {\sl flip},
the action of lifting two dominoes forming a square 
or three lozenges forming a hexagon
and placing them back after a rotation of $90^\circ$ or $60^\circ$.
Two tilings are joined by a flip
iff their height functions are different in one vertex only;
this difference must be 3 (resp. 4) 
in the triangulated (resp. quadriculated) case.
Equivalently, in the pile-of-cubes (or pile-of-blocks) interpretation,
a flip corresponds to adding or removing a cube (or block).
This suggests a basic property of $T$
when $R$ is simply connected:

{\nobf Theorem 3.1: }
{\sl If $R$ is a simply connected region,
$T$ is connected.}

The inductive step is to prove that
if a tiling $t_1$ is larger than a tiling $t_2$,
we may perform a flip on $t_1$ to obtain $t_3$
with $t_1 > t_3 \ge t_2$.
Let $\theta_i$ be the height function for $t_i$;
remember that all height functions coincide along the boundary.
To find out where to flip,
consider the set of vertices where $\theta_1 - \theta_2$
is maximal and, within this set,
choose the vertex with maximum $\theta_1$.
Since all height functions coincide in the boundary,
this point is a local maximum of $\theta_1$.
By the characterization of height functions,
a flip is possible at that point
(we found a removable cube or block).

The triangulated version of theorem 2.2 follows from theorem 3.1:
for triangulated surfaces, flips preserve parity of the monomials
in the determinant expansion of $B$.

As usual, define distance in the graph $T$
as the minimum number of edges (flips) in a path joining two vertices.
The following distance formula follows from the same
inductive step.

{\nobf Theorem 3.2: }
{\sl 
Let $R$ be a simply connected region.
The distance between tilings $t_1$ and $t_2$,
with height functions $\theta_1$ and $\theta_2$, is
$$d(t_1,t_2) = {1 \over X}\sum_v{|\theta_1(v) - \theta_2(v)|},$$
where $X=3$ or 4 depending whether $R$
is triangulated or quadriculated.}

As usual, we call such minimal paths {\sl geodesics}.
Some geodesics are pretty much the same:
we may freely change the order of two successive flips
($f_1$ and $f_2$)
if they act on disjoint hexagons or squares.
We call two such geodesics equivalent
and glue a 2-cell (along the loop $f_1f_2f_1^{-1}f_2^{-1}$)
to the graph $T$ to make them homotopic.
Similarly, glue $k$-cells to identify
all permutations of $k$ independent flips,
obtaining from $T$ the $CW$-complex $\cal T$.

{\nobf Theorem 3.3: }
{\sl If $R$ is simply connected, $\cal T$ is contractible.}

As a consequence, all geodesics are equivalent
and all closed curves null-homotopic.

We now consider non-simply connected planar regions $R$.
The graph $T$ is now usually disconnected:
to see this, we consider the following invariant under flips.
Take a {\sl cut}, i.e., a simple oriented curve in $R$ consisting 
of edges with extrema in different boundary components of $R$.
If $R$ is triangulated (resp. quadriculated),
define the {\sl flow} of a tiling across a cut
as the number of lozenges (resp. dominoes) crossing the cut,
counted with a sign determined by the colour
of the triangle (resp. square) to the left of the cut.
For instance, the two tilings in figure 3.1
are not in the same connected component of $T$
because the flows across the indicated cut are 1 and 0.
The reader can easily see that homologous cuts
yield equivalent conditions,
the flows differing by a fixed constant from one cut to the other.
In terms of height sections,
it turns out that the height bundle does not depend on the tiling
and two tilings have the same flow across all cuts
iff their height sections coincide on the whole boundary.

\midinsert
{
\smallskip

\centerline{\psfig{file=fig31.eps,width=3.5in}}


\smallskip
\centerline{\eightbf Figure 3.1}

}
\endinsert

{\nobf Theorem 3.4: }
{\sl
Let $R$ be a quadriculated planar region.
Two tilings of $R$ belong to the same connected component of $T$
iff their flows coincide across all cuts.
For two tilings in the same connected component of $T$,
the distance formula in theorem 3.2 holds.
Furthermore, connected components of $\cal T$
(constructed as before) are contractible.}

It would be interesting to see what happens
for triangulated planar regions.

A further generalization consists in considering
arbitrary triangulated or quadriculated surfaces.
The cut and flow invariants still make sense
if the surface is orientable and bicoloured
(no longer necessary conditions!).
The appropriate generalization to the general case
involves homology and cohomology with local coefficients,
which we shall not discuss.
In figure 3.2, we show that cut and flow conditions
are no longer sufficient,
even assuming orientability and bicolourability.
For obvious reasons, we call these configurations {\sl ladders};
rather surprisingly, their exact position turns out
to be the only additional condition for domino tilings
to be in the same connected component of $T$.
As a further complication, if the surface is a (bicolourable) torus,
certain connected components of $\cal T$
are no longer contractible,
but homotopically equivalent to $\SS^1$
and therefore contain ``closed geodesics''.

\midinsert
{
\smallskip

\centerline{\psfig{file=fig32.eps,width=3.5in}}


\smallskip
\centerline{\eightbf Figure 3.2}

}
\endinsert

{\nobf Theorem 3.5: }
{\sl
Let $R$ be an arbitrary quadriculated surface.
If the boundary of $R$ is non-empty,
connected components of $\cal T$ are contractible.
In general, they are either contractible
or homotopically equivalent to circles.
}

\goodbreak
\bigskip
{\nobigbf References}

\nobreak
\smallskip
\nobreak

{\parindent=40pt

\item{[B]}{R. Berger,
{The undecidability of the domino problem,}
Mem. Amer. Math. Soc. {\bf 66} (1966).}

\item{[BK]}{E. A. Bender and D. E. Knuth,
{Enumeration of plane partitions,}
J. Comb. Theor. {\bf A13}, 40-54 (1972).}

\item{[CL]}{J. H. Conway and J. C. Lagarias,
{Tilings with polyominoes and combinatorial group theory,}
J. Comb. Theor. {\bf A53}, 183-208 (1990).}

\item{[DaT]}{G. David and C. Tomei,
{The problem of the calissons,}
Amer. Math. Monthly, {\bf 96}, 429-431 (1989).}

\item{[DeT]}{P. A. Deift and C. Tomei,
{On the determinant of the adjacency matrix for a planar
sublattice,} J. Comb. Theor. {\bf B35}, 278-289 (1983).}

\item{[GJ]}{M. Garey and D. S. Johnson,
{Computers and intractability:
a guide to the theory of NP-completeness,}
Freeman, San Francisco, 1979.}

\item{[K]}{P. W. Kasteleyn,
{The statistics of dimers on a lattice I. The number of dimer
arrangements on a quadratic lattice,}
Phisica {\bf 27}, 1209-1225 (1961).}

\item{[LL]}{E. H. Lieb and M. Loss,
{Fluxes, Laplacians and Kasteleyn's theorem,}
Duke Math. Jour., {\bf 71}, 337-363 (1993).}

\item{[M]}{P. A. MacMahon,
{Combinatory analysis,} vol. 2, Cambridge University Press, 1916;
reprinted by Chelsea, New York, 1960.}

\item{[STCR]}{N. C. Saldanha, C. Tomei,
M. A. Casarin Jr. and D. Romualdo,
{Spaces of domino tilings,}
Discr. Comp. Geom., to appear.}

\item{[T]}{W. P. Thurston,
{Conway's tiling groups,}
Amer. Math. Monthly, {\bf 97}, 8, 757-773 (1990).}

}

\bigskip\bigskip

\goodbreak

\line{Nicolau C. Saldanha, IMPA and PUC-Rio \hfill}
\line{e-mail: nicolau@impa.br \hfill}

\smallskip

\line{Carlos Tomei, IMPA and PUC-Rio \hfill}
\line{e-mail: tomei@impa.br \hfill}

\medskip

\line{Instituto de Matematica Pura e Aplicada,
Estrada Dona Castorina, 110 \hfill}
\line{Jardim Bot\^anico, Rio de Janeiro, RJ 22460-320, BRAZIL\hfill}

\smallskip

\line{Departamento de Matem\'atica, PUC-Rio,
Rua Marqu\^es de S\~ao Vicente, 225 \hfill}
\line{G\'avea, Rio de Janeiro, RJ 22453-900, BRAZIL \hfill}

\bye